\documentclass[a4paper,12pt]{article}

\usepackage[english]{babel}

\usepackage[utf8]{inputenc}

\usepackage{fancyhdr}
\usepackage{amsthm}
\theoremstyle{definition}
\newtheorem{defn}{Definition}[section]
\newtheorem{thm}[defn]{Theorem}
\newtheorem{cor}[defn]{Corollary}
\newtheorem{lem}[defn]{Lemma}
\newtheorem{oss}[defn]{Remark}
\newtheorem*{quest}{Question}

\usepackage{amsmath}

\usepackage{amssymb}

\usepackage{amsbsy}

\usepackage{bm}

\usepackage{float}

\usepackage[usenames,dvipsnames]{xcolor}

\usepackage[nottoc]{tocbibind}

\usepackage[bookmarks,bookmarksopen]{hyperref}
\hypersetup{pdfauthor={Raffaella Cutolo},pdftitle={Choiceless large cardinals and set-theoretic potentialism}}
\hypersetup{citebordercolor=ForestGreen}
\hypersetup{linkbordercolor=RoyalBlue}

\begin{document}

\centerline{\Large{\bf Choiceless large cardinals and}}
\medskip
\centerline{\Large{\bf set-theoretic potentialism}}

\bigskip
\centerline{Raffaella Cutolo and Joel David Hamkins}

\medskip
\centerline{\today}

\medskip
\begin{abstract}
We define a \emph{potentialist system} of ZF-structures, that is, a collection of \emph{possible worlds} in the language of ZF connected by a binary \emph{accessibility relation}, achieving a so called ``potentialist account" of the full background set-theoretic universe $V$. The definition involves Berkeley cardinals, the strongest known large cardinal axioms, inconsistent with the Axiom of Choice. In fact, as background theory we assume just ZF. It turns out that the propositional modal assertions which are valid at every world of our system are exactly those in the modal theory S4.2. Moreover, we characterize the worlds satisfying the potentialist maximality principle, and thus the modal theory S5, both for assertions in the language of ZF and for assertions in the full potentialist language.
\end{abstract}

\section{Introduction}

In the current scenario of set theory, we are faced with a conflict between large cardinal axioms and the Axiom of Choice. In fact, there is a whole new hierarchy of ZF large cardinals - the Berkeley hierarchy - which contradict AC and lie beyond the Kunen inconsistency of Reinhardt cardinals. Such ``choiceless" large cardinals have been recently introduced in \cite{BKW} and the investigation of their consistency is very involved in the present main foundational questions concerning the universe of set theory. But let us point out something else that is of interest here, namely: if we drop AC then the set-theoretic universe $V$ grows upward. This observation raises a \emph{potentialist} perspective, that is, one in which the universe of set theory reveals gradually, and never completely, as we progressively take under consideration new fragments of it; indeed, we can actually think of we access higher and higher parts of the set-theoretic universe by considering stronger and stronger large cardinals.

Recent works of the second author focus on the idea of set-theoretic potentialism and the analysis of the modal principles validated by specific potentialist systems. The general definition is stated below.

\begin{defn}
A $\bf potentialist\ system$ is a collection $\mathcal{W}$ of structures in a common language $\mathcal{L}$ called ``worlds", equipped with a binary accessibility relation $\mathcal{R}$, such that:
\begin{itemize}
\item $\mathcal{R}$ is reflexive and transitive;
\item whenever $M\mathcal{R}N$, then $M$ is - or embeds to - a substructure of $N$.
\end{itemize}
\end{defn}

\noindent So, a potentialist system is a Kripke model of $\mathcal{L}$-structures for some language $\mathcal{L}$. In order to study how truth of an assertion $\varphi$ propagates through the worlds of $\mathcal{W}$, one adds to the basic language $\mathcal{L}$ the modal operators $\Diamond$ and $\Box$, expressing, respectively, the notions of \emph{possibility} and \emph{necessity}:
\begin{itemize}
\item $\Diamond\varphi$ holds at a world $M$ (that is, ``$\varphi$ is possible over $M$") if $\varphi$ holds at some world $N$ such that $M\mathcal{R}N$;
\item $\Box\varphi$ holds at a world $M$ (that is, ``$\varphi$ is necessary over $M$") if $\varphi$ holds at all worlds $N$ such that $M\mathcal{R}N$.
\end{itemize}
Now one can ask which propositional modal assertions are \emph{valid} in the whole system $\mathcal{W}$ (that is, hold in every world of $\mathcal{W}$); the point is that determining the modal validities of a potentialist system gives a precise account of how its worlds interact with respect to their respective truths.

Let us turn to our particular case, whose hallmark is to combine choiceless large cardinals with the potentialist ideas. Indeed, we consider the concept of set-theoretic potentialism that arises from elementary embeddings of a transitive set into itself, where we view $M$ as \emph{accessing} $N\supseteq M$ whenever the restriction to $M$ of any elementary embedding $j:N\to N$ yields an elementary embedding $j':M\to M$. Such a definition of the accessibility relation results in an interesting case as in the context of Berkeley cardinals, one can arrange \emph{non-trivial} elementary embeddings fixing any desired set. The key point is that every given set is definable in some big transitive set and if there is a Berkeley cardinal $\delta$ then, by definition, any transitive set $M$ containing $\delta$ as a member admits non-trivial elementary embeddings $j:M\to M$, whose critical points are in fact cofinal in $\delta$.

\begin{defn}
A cardinal $\delta$ is a $\bf Berkeley\ cardinal$ if for every transitive set $M$ such that $\delta\in M$, and for every ordinal $\eta<\delta$, there exists a non-trivial elementary embedding $j:M\to M$ with $\eta<crit(j)<\delta$.\footnote{In the choiceless context, being non-trivial means $j$ is not the identity on the ordinals.}
\end{defn}

It turns out that the set-theoretic universe $V$ equals the union of the worlds of our potentialist system, and given any world $M$ and any set $a$, there is a world $N$ \emph{accessed} by $M$ such that $a\in N$. Thus, truth in $V$ is approximated by truth in our worlds: we can assert any property concerning any set of $V$ from any of the worlds of our system by using the diamond operator, and we can progressively move from any world to the wider perspective of another world which is ``closer" to $V$ in that it contains additional sets and is capable to satisfy additional properties about them. The primary goal here will be to provide a definite account and determine the valid modal principles of this kind of set-theoretic potentialism; but let us mention that as a further perspective, maybe one could use such a multiverse setting to investigate further the fundamental question of the consistency of the choiceless large cardinals.

For the basics of set-theoretic potentialism and the various potentialist systems analyzed so far refer to \cite{Ham}. For more on the choiceless large cardinals see \cite{BKW}.

\section{A potentialist account of the set-theoretic universe \!\textit{V}}

We start with a preliminary lemma motivating the definition of the accessibility relation we shall consider.

\begin{lem}\label{liftlemma}
For every transitive set $M$, for every set $a$, there exists a transitive set $N\supseteq M$ with $a\in N$ such that every elementary embedding $j:N\to N$ lifts some elementary embedding $j':M\to M$.
\end{lem}

\proof
Let $M$ be a transitive set and let $a$ be any set. As shown in \cite{BKW}, there exists a transitive set $N$ such that $M,\,a\in N$ and $M$ is definable (without parameters) in $N$. Thus, $M\subseteq N$ and every $j:N\to N$ fixes $M$, which implies $j``M\subseteq M$. Therefore $j\restriction M:M\to M$, and so actually $j$ lifts $j'=j\restriction M:M\to M$.
\endproof

\begin{defn}
Let $\delta$ be a Berkeley cardinal.
\begin{itemize}
\item Let $\mathcal{M}_\delta=\left\{M:M\ is\ transitive\wedge\delta\in M\right\}$.
\item Let $\mathcal{R}$ be the binary relation on $\mathcal{M}_\delta$ defined as follows: for $M,\,N\in\mathcal{M}_\delta$,
\begin{center}
$M\mathcal{R}N$ iff $M\subseteq N$ and every elementary embedding $j:N\to N$ lifts some elementary embedding $j':M\to M$; that is, $M\mathcal{R}N$ iff $M\subseteq N$ and for every elementary embedding $j:N\to N$, $j\restriction M:M\to M$.
\end{center}
\end{itemize}
\end{defn}

\noindent Note that the assumption that there is a Berkeley cardinal $\delta$ and the choice of $\mathcal{M}_\delta$ as collection of worlds ensure that every world $M$ admits non-trivial elementary embeddings $j:M\to M$, so $\mathcal{R}$ is not merely reduced to the subset relation. It is trivial that $\mathcal{R}$ is reflexive; also, $\mathcal{R}$ is transitive: in fact, if $j:H\to H$ lifts $j':N\to N$, which in turn lifts $j'':M\to M$, then $j$ lifts $j''$, as $j\restriction M=(j\restriction N)\restriction M=j'\restriction M=j''$. Therefore, $\langle\mathcal{M}_\delta,\mathcal{R}\rangle$ is a \emph{potentialist system} of ZF-structures. Since $V_\alpha\in\mathcal{M}_\delta$ for any $\alpha>\delta$, we have that $V=\bigcup_{M\in\mathcal{M}_\delta}M$. Moreover, we show that $\mathcal{M}_\delta$ provides a \emph{potentialist account} of the set-theoretic universe $V$, meaning that every world in $\mathcal{M}_\delta$ is a substructure of $V$ and for every $M\in\mathcal{M}_\delta$ and every set $a$ there is a world $N\in\mathcal{M}_\delta$ accessed by $M$ such that $a\in N$.

\begin{lem}\label{accountV}
$\mathcal{M}_\delta$ provides a potentialist account of the universe $V$.
\end{lem}

\proof
First, for every $M\in\mathcal{M}_\delta$, $\langle M,\in\rangle$ is a substructure of $\langle V,\in\rangle$ (i.e., $M\subseteq V$ and $\in^M=\,\in^V\restriction M$). Further, if $M\in\mathcal{M}_\delta$ and $a\in V$, then by Lemma \ref{liftlemma} there exists a transitive set $N$ such that $\left\{M,a\right\}\in N$ and for every elementary embedding $j:N\to N$, $j\restriction M:M\to M$. Since $\delta\in M\subseteq N$, we have that $N\in\mathcal{M}_\delta$; since $a\in N$ and $N$ is accessed by $M$, we are done.
\endproof

\begin{oss}
Notice that:
\begin{enumerate}
\item By Lemma \ref{liftlemma}, for every $M\in\mathcal{M}_\delta$ there exist cofinally many $N\in\mathcal{M}_\delta$ which are accessed by $M$, meaning that such $N$ can accommodate any given set.
\item In particular, for every $M\in\mathcal{M}_\delta$, for every set $a$, there exists a rank initial segment $V_\alpha\supseteq M$ with $a\in V_\alpha$ such that every elementary embedding $j:V_\alpha\to V_\alpha$ lifts some elementary embedding $j':M\to M$, that is, such that $V_\alpha$ is accessed by $M$.
\end{enumerate}
\end{oss}

\noindent For any assertion $\varphi$ in the language of ZF, the \emph{potentialist translation} $\varphi^\diamond$ is the assertion in the potentialist language $\mbox{ZF}^\diamond$ (which augments the language of ZF with the modal operators $\Diamond$ and $\Box$) achieved by replacing every instance of $\exists x$ with $\Diamond\exists x$ and every instance of $\forall x$ with $\Box\forall x$. As an immediate corollary of Lemma \ref{accountV}, we get that truth in $V$ is equivalent to potentialist truth at the worlds of $\mathcal{M}_\delta$.

\begin{cor}
For any ZF-formula $\varphi$ and for any $a_0,\dots,a_n\in V$, we have:
\begin{center}
$V\models\varphi(a_0,\dots,a_n)$ iff $M\models_{\mathcal{M}_\delta}\varphi^\diamond(a_0,\dots,a_n)$,
\end{center}
for any $M\in\mathcal{M}_\delta$ in which $a_0,\dots,a_n$ exist.
\end{cor}

\noindent Let us now state what it precisely means for a modal assertion to be \emph{valid} with respect to our potentialist system.

\begin{defn}
A modal assertion $\varphi(p_0,\dots,p_n)$ in the language of propositional modal logic is $\bf valid$ at a world $M$ in $\mathcal{M}_\delta$ for a certain class of assertions, if all the resulting substitution instances $\varphi(\psi_0,\dots,\psi_n)$, where assertion $\psi_i$ from the allowed class is substituted for the propositional variable $p_i$, are true at $M$. $\varphi$ is valid in $\mathcal{M}_\delta$ if it is valid at every world of $\mathcal{M}_\delta$.
\end{defn}

\noindent Of course, the main question arising here is the following:

\begin{quest}
What is the modal logic of $\mathcal{M}_\delta$? That is, which are the modal principles valid in $\mathcal{M}_\delta$?
\end{quest}

\section{The modal logic of \texorpdfstring{$\mathcal{M}_\delta$}{M_delta}}

In this section we provide lower and upper bounds on the modal validities of $\mathcal{M}_\delta$, and finally prove that the modal logic of $\mathcal{M}_\delta$ is exactly S4.2.

\begin{defn}
The modal theory S4 is obtained from the following axioms by closing under modus ponens and necessitation.
\begin{itemize}
\item (K) $\Box(\varphi\to\psi)\to (\Box\varphi\to\Box\psi)$
\item (Dual) $\neg\Diamond\varphi\leftrightarrow\Box\neg\varphi$
\item (S) $\Box\varphi\to\varphi$
\item (4) $\Box\varphi\to\Box\Box\varphi$
\end{itemize}
\end{defn}

\begin{thm}
The modal theory S4 is valid at every world of $\mathcal{M}_\delta$.\footnote{In fact, every potentialist system validates S4.}\footnote{Unless otherwise specified, the validities hold for all assertions in $\mbox{ZF}^\diamond$, with parameters.}
\end{thm}

\proof
Let $M\in\mathcal{M}_\delta$.
\begin{itemize}
\item (K). Suppose $\Box(\varphi\to\psi)$ and $\Box\varphi$ hold in $M$. Then, $\varphi\to\psi$ and $\varphi$ hold in any world $N$ accessed by $M$. Therefore, by modus ponens, $\psi$ holds in any such $N$, that is, $\Box\psi$ holds in $M$.
\item (Dual). Immediate.
\item (S). Follows immediately from the fact that every world accesses itself.
\item (4). If $\varphi$ holds in any world $N$ accessed by $M$, then so does $\Box\varphi$, as any world $H$ accessed by $N$ is also accessed by $M$.
\end{itemize}
\endproof

\begin{defn}
The modal theory S4.2 is obtained from S4 by adding the axiom (.2) $\Diamond\Box\varphi\to\Box\Diamond\varphi$.
\end{defn}

\begin{thm}\label{S4.2}
The modal theory S4.2 is valid at every world of $\mathcal{M}_\delta$.
\end{thm}

\proof
(.2). Let $M\in\mathcal{M}_\delta$. Assume $\Diamond\Box\varphi$ holds in $M$, that is, there exists $N\in\mathcal{M}_\delta$ such that $M\mathcal{R}N$ and $\Box\varphi$ holds in $N$. Let $H\in\mathcal{M}_\delta$ be such that $M\mathcal{R}H$. We need to show that $\Diamond\varphi$ holds in $H$. Note that there exists a transitive set $K$ such that $\langle N,H\rangle$ is definable in $K$. Since $K$ is transitive and $N,\,H\subseteq K$, we have that $\delta\in K$ and so $K\in\mathcal{M}_\delta$. Now, take any non-trivial elementary embedding $j:K\to K$. Then, $j(N)=N$ and $j(H)=H$. So, $j\restriction N:N\to N$ and $j\restriction H:H\to H$. Therefore, $N$ and $H$ both access $K$. Since $K$ is accessed by $N$, $K$ satisfies $\varphi$; but then, since $K$ is accessed by $H$ and $\varphi$ holds in $K$, $\Diamond\varphi$ holds in $H$.
\endproof

\noindent Theorem \ref{S4.2}, whose proof actually shows our potentialist system is ``convergent" - or ``locally directed"\footnote{In fact, it shows that whenever $M\mathcal{R}N$ and $M\mathcal{R}H$ then there exists $K$ such that $N\mathcal{R}K$ and $H\mathcal{R}K$, and so, that $\mathcal{M}_\delta$ has \emph{amalgamation}.} - and therefore validates (.2), establishes a first significant lower-bound result. In order to provide upper bounds on the validities of $\mathcal{M}_\delta$, and then determine the exact set of modal principles valid through the whole system, we recall the definitions of \emph{switches}, \emph{buttons} and \emph{dials}, specific kinds of control statements first introduced in \cite{HL}; in particular, we will be interested in finding assertions satisfying such definitions which also have the property of being \emph{independent}.

\begin{defn}
An assertion $s$ is a $\bf switch$ if both $\Diamond s$ and $\Diamond\neg s$ are true at every world, that is, both $\Box\Diamond s$ and $\Box\Diamond\neg s$ hold. $s$ is a switch at a particular world $M$ if $\Diamond s$ and $\Diamond\neg s$ are true at all the worlds accessed by $M$.
\end{defn}

\begin{defn}
A $\bf button$ is a statement $b$ such that $\Diamond\Box b$ is true at every world, that is, $\Box\Diamond\Box b$ holds. The button is \emph{pushed} at a world if $\Box b$ holds at that world, and otherwise unpushed.
\end{defn}

\begin{defn}
A (possibly infinite) list of statements $d_0,d_1,d_2,\dots$ is a $\bf dial$ if every world satisfies exactly one of the statements $d_i$ and every world can access another world with any prescribed dial value. If a world satisfies $d_i$, then we say that the dial value is $i$ in that world.
\end{defn}

\begin{defn}
A family of switches is $\bf independent$ if one can always flip the truth values of any finitely many of the switches so as to realize any desired finite pattern of truth.
\end{defn}

\begin{defn}
A family of buttons and switches is $\bf independent$ if there is a world at which the buttons are unpushed, and every world $M$ accesses a world $N$ in which any additional button may be pushed without pushing any other as-yet unpushed button from the family, while also setting any finitely many of the switches so as to have any desired pattern in $N$; and similarly with dials.
\end{defn}

\begin{defn}
The modal theory S5 is obtained from S4 by adding the axiom (5) $\Diamond\Box\varphi\to\varphi$, which we call \emph{potentialist maximality principle} (MP).
\end{defn}

\noindent The following theorem summarizes some key results - first proved in \cite{HL}, and developed further in \cite{HLL} - we shall use.

\begin{thm}
The following hold.
\begin{enumerate}
\item If $\mathcal{W}$ is a Kripke model and a world $M\in\mathcal{W}$ admits arbitrarily large finite collections of independent switches, then the propositional modal assertions valid at $M$ are contained in the modal theory S5. In particular, if the switches work throughout $\mathcal{W}$, then the validities of every world of $\mathcal{W}$, and so the validities of $\mathcal{W}$, are contained within S5.
\item A Kripke model $\mathcal{W}$ admits arbitrarily large finite families of independent switches if and only if it admits arbitrarily large finite dials.
\item If $\mathcal{W}$ is a Kripke model that admits arbitrarily large finite families of independent buttons and switches, or independent buttons independent of a dial, then the propositional modal validities of $\mathcal{W}$ are contained in S4.2. The validities of any particular world in which the buttons are not yet pushed are contained in S4.2, and in any case, are contained in S5.
\end{enumerate}
\end{thm}

\begin{thm}\label{S5}
The propositional modal validities of $\mathcal{M}_\delta$ are contained in the modal theory S5.
\end{thm}

\proof
It suffices to show that $\mathcal{M}_\delta$ admits arbitrarily large finite dials. We shall show that it admits in fact an infinite dial (notice that from an infinite dial, we can construct finite dials of any given size by keeping any desired finite number of dial statements and adding the statement that none of them holds). For $i<\omega$, let $d_i$ be the assertion that the height of the ordinals is $\lambda+i$, where $\lambda$ is a limit ordinal or zero. These statements are expressible in the language of ZF (without parameters or modal vocabulary), correctly interpreted inside any transitive set, and so inside any world $M\in\mathcal{M}_\delta$. Let us show that they form a dial. First, since any ordinal is uniquely expressed as $\lambda+i$ for some limit ordinal $\lambda$ or zero and some finite $i<\omega$, every world in $\mathcal{M}_\delta$ satisfies exactly one of the statements $d_i$. It remains to prove that every world can access another world with any desired dial value. Let $M\in\mathcal{M}_\delta$ and fix $i<\omega$. Let $V_\theta$ be such that $M\subseteq V_\theta$. Let $N=V_\theta\cup\mbox{tr\,cl}(\left\{\langle\theta,M\rangle\right\})\cup(\theta+i)$. Then $N\in\mathcal{M}_\delta$ and $M$ is definable in $N$, which implies $M\mathcal{R}N$, and the dial value in $N$ is $i$.
\endproof

\begin{thm}
$\mathcal{M}_\delta$ satisfies exactly S4.2, that is, the modal logic of $\mathcal{M}_\delta$ is S4.2.
\end{thm}

\proof
We show that $\mathcal{M}_\delta$ admits arbitrarily large finite families of independent buttons independent of a dial, which implies that the modal validities of $\mathcal{M}_\delta$ are contained in, and hence by Theorem \ref{S4.2} equal to, S4.2. For $i<\omega$, let $d_i$ be the assertion that the height of the ordinals is $\lambda+i$, where $\lambda$ is a limit ordinal or zero; we already showed that these statements form a dial. For $m<\omega$, let $b_m$ be the assertion that the set $m\cdot\mathbb{N}=\left\{m\cdot k:k<\omega\right\}$ exists; let us show that these statements are independent buttons independent of the above dial. Since for every $m<\omega$, if the assertion $b_m$ is true in some transitive set then it will continue to be true in any larger transitive set, each $b_m$ is a button. These buttons are independent because every world $M$ accesses a world $N$ in which any additional button $b_m$ may be pushed without pushing any other as-yet unpushed button from the family. Finally, the above buttons and dial values can be controlled independently of each other.
\endproof

\noindent By Theorem \ref{S5}, S5 is a definite upper bound on the validities of $\mathcal{M}_\delta$. An interesting point would therefore be to determine which worlds of $\mathcal{M}_\delta$ realize the maximum set of validities. In other words:

\begin{quest}
Which worlds of $\mathcal{M}_\delta$ satisfy the potentialist maximality principle?
\end{quest}

\section{The worlds satisfying MP}

We now give a characterization of the worlds of $\mathcal{M}_\delta$ satisfying S5. Depending on the language we consider, we get different criteria. The following concepts are involved.

\begin{defn}
An ordinal $\theta$ is $\boldsymbol{\Sigma_n}$-$\bf correct$ if $V_\theta\prec_{\Sigma_n}V$, meaning that $V_\theta$ and $V$ agree on the truth of $\Sigma_n$ formulas with parameters from $V_\theta$.
\end{defn}

\begin{defn}
A cardinal $\theta$ is $\bf correct$ if it is $\Sigma_n$-correct for every $n$, that is, if it realizes the scheme $V_\theta\prec V$.\footnote{Note: this concept is not expressible as a single assertion in the language of set theory, although it can be expressed as a scheme of statements.}
\end{defn}

\begin{thm}\label{sigma2}
The following are equivalent.
\begin{enumerate}
\item The potentialist maximality principle holds in a world $M\in\mathcal{M}_\delta$ for assertions in the language of ZF with parameters from $M$.
\item $M=V_\theta$ for some $\Sigma_2$-correct cardinal $\theta>\delta$.
\end{enumerate}
\end{thm}

\proof
Let $M\in\mathcal{M}_\delta$.
\begin{itemize}
\item ($1\Rightarrow 2$). Assume $M$ satisfies (5) $\Diamond\Box\varphi\to\varphi$ for assertions in the language of ZF with parameters from $M$. First, note that $M$ must be a $V_\theta$ with $\theta$ limit. In fact, for all $a\in M$, the existence of the power set of $a$ is possibly necessary, and it follows that $M$ thinks $\mathcal{P}(a)$ exists, and that $M$ computes the power sets correctly; moreover, $M$ computes $V_\alpha$ correctly for any ordinal $\alpha\in M$, since the existence of $V_\alpha$ is possibly necessary. Also, for any set $a$, it is possibly necessary that $a\in V_\alpha$ for some ordinal $\alpha$, and so this is already true in $M$. Now, suppose $\varphi(\vec{a})$ is a $\Sigma_2$ statement true in $V$, with $\vec{a}\in M$. This is witnessed by the existence of some ordinal $\alpha$ for which $V_\alpha$ satisfies $\psi(\vec{a})$ for some assertion $\psi$. So, it is possibly necessary the statement that there is an ordinal $\alpha$ for which $V_\alpha$ exists and satisfies $\psi(\vec{a})$. Thus, by (5), this statement must be true in $M$, and so $\varphi(\vec{a})$ is true in $M$. Therefore, $\theta$ is $\Sigma_2$-correct.
\item ($2\Rightarrow 1$). Assume $M=V_\theta$ for $\theta>\delta$ a $\Sigma_2$-correct cardinal. Suppose $M$ satisfies $\Diamond\Box\varphi(\vec{a})$, $\vec{a}\in M$. Then, there exists $N\in\mathcal{M}_\delta$ accessed by $M$ that satisfies $\Box\varphi(\vec{a})$. Since $\Diamond\Box\varphi(\vec{a})$ is a $\Sigma_2$ statement in $V$, it must be true inside $M$, being $M$ $\Sigma_2$-correct. So, there exists $m\in M$ such that $m$ satisfies $\Box\varphi(\vec{a})$. Without loss of generality, there exists $m=V_\alpha$ like this (inside $M$), and so the smallest one is definable and so it accesses $M$. Since the statement that $V_\alpha$ satisfies $\Box\varphi(\vec{a})$ is a $\Pi_1$ statement, it is absolute between $M$ and $V$ by the $\Sigma_2$-correctness of $M$, and so, it holds in $V$. Thus, $\varphi(\vec{a})$ holds in $M$, which therefore satisfies (5).
\end{itemize}
\endproof

\begin{thm}\label{MP}
The following are equivalent.
\begin{enumerate}
\item The potentialist maximality principle holds in a world $M\in\mathcal{M}_\delta$ for assertions in the potentialist language $\mbox{ZF}^\diamond$ with parameters from $M$.
\item $M=V_\theta$ for some correct cardinal $\theta>\delta$.
\end{enumerate}
\end{thm}

\proof
Let $M\in\mathcal{M}_\delta$.
\begin{itemize}
\item ($1\Rightarrow 2$). Assume $M$ satisfies (5) $\Diamond\Box\varphi\to\varphi$ for assertions in the potentialist language $\mbox{ZF}^\diamond$ with parameters from $M$. Then by Theorem \ref{sigma2}, $M=V_\theta$ for $\theta$ $\Sigma_2$-correct. By the potentialist translation, truth in $V$ is expressible as $\varphi^\diamond$-truth in $M$. Thus, $M$ can express the statement that there exists a $\Sigma_n$-correct cardinal (as this is definable in $V$). For each $n$, this statement is a button in $\mbox{ZF}^\diamond$. So $M$ is a limit of $\Sigma_n$-correct cardinals, and therefore $M$ is fully correct.
\item ($2\Rightarrow 1$). Assume $M=V_\theta$ for $\theta>\delta$ a correct cardinal. Suppose $M$ satisfies $\Diamond\Box\varphi(\vec{a})$, where $\vec{a}\in M$ and $\varphi$ is a $\mbox{ZF}^\diamond$ assertion. Then, there exists $N\in\mathcal{M}_\delta$ accessed by $M$ that satisfies $\Box\varphi(\vec{a})$. Since the existence of such a set $N$ and the potentialist semantics are expressible in the language of set theory, it follows from $M\prec V$ that there is such a set inside $M$. So, there exists $m\in M$ such that $m$ satisfies $\Box\varphi(\vec{a})$. Without loss of generality, this $m$ has form $V_\alpha$, and so the smallest one is definable and therefore accesses $M$. Thus, $\varphi(\vec{a})$ holds in $M$.
\end{itemize}
\endproof

\begin{oss}
Note that one can view Theorem \ref{MP} as a ZF theorem scheme asserting the equivalence of two schemes; or, one could view it as a theorem of G\"odel-Bernays set theory augmented with the assumption that there is a predicate for first-order set-theoretic truth (that theory is provable, for example, in Kelley-Morse set theory).
\end{oss}

\section{Consistency of MP and conclusive remarks}

Theorem \ref{sigma2} and Theorem \ref{MP} characterize the worlds of $\mathcal{M}_\delta$ satisfying MP, respectively, for assertions in the language of set theory and for assertions in the full potentialist language, with parameters. Let us remark that there exist indeed such worlds in $\mathcal{M}_\delta$. In fact, for the first case, note that by the L\'evy-Montague reflection theorem (which is a ZF result), the class of all $\Sigma_2$-correct cardinals is closed and unbounded in the ordinals. For the second case, observe that although the existence of a correct cardinal is not provable in ZF, it is relatively consistent with ZF (see \cite{MP}); so, it is relatively consistent with ZF that there exist worlds in $\mathcal{M}_\delta$ satisfying MP for assertions in the potentialist language $\mbox{ZF}^\diamond$.

Recall that the definition of our potentialist system leverages on the assumption that there is a Berkeley cardinal $\delta$; one may ask if there is any world in $\mathcal{M}_\delta$ which satisfies that there exists a Berkeley cardinal, and the answer is yes: in fact, the property of being a Berkeley cardinal is $\Pi_2$, so for any $\Sigma_2$-correct cardinal $\theta$, $V_\theta$ correctly recognizes the Berkeley cardinals below $\theta$. In other words, inside all the worlds of $\mathcal{M}_\delta$ satisfying MP, $\delta$ itself is still a Berkeley cardinal; but these are not the only worlds in $\mathcal{M}_\delta$ recognizing $\delta$ is Berkeley: in fact, as noted in \cite{BKW}, if $\delta$ is a Berkeley cardinal then for all limit ordinals $\lambda>\delta$, $V_\lambda$ thinks that $\delta$ is Berkeley.

Finally, since every world $M$ in $\mathcal{M}_\delta$ can access a $V_\lambda$ with $\lambda$ limit, the assertion $\varphi_{\scriptsize{\mbox{BC}}}$ that there exists a Berkeley cardinal is possible over any $M$, that is, $\Diamond\varphi_{\scriptsize{\mbox{BC}}}$ holds at every world.

\end{document}